\documentclass[11pt]{article}
\usepackage{geometry}                
\usepackage{graphicx}
\usepackage{amssymb}
\usepackage{epstopdf}
\usepackage{amsmath,amsthm,amscd,amssymb}
\usepackage{latexsym}
\usepackage[colorlinks,citecolor=red,pagebackref,hypertexnames=false]{hyperref}

\numberwithin{equation}{section}

\theoremstyle{plain}
\newtheorem{theorem}{Theorem}[section]
\newtheorem{lemma}[theorem]{Lemma}

\theoremstyle{definition}
\newtheorem{definition}[theorem]{Definition}

\theoremstyle{remark}
\newtheorem{remark}[theorem]{Remark}

\newtheorem{case[theorem]}{Case}

\title{Generalized incidence theorems, homogeneous forms and sum-product estimates in finite fields}
\author{David Covert, Derrick Hart, Alex Iosevich, Doowon Koh, and Misha Rudnev}

\begin{document}
\maketitle

\begin{abstract} In recent years, sum-product estimates in Euclidean space and finite fields have been studied using a variety of combinatorial, number theoretic and analytic methods. Erdos type problems involving the distribution of distances, areas and volumes have also received much attention. In this paper we prove a relatively straightforward function version of an incidence results for points and planes previously established in \cite{HI07} and \cite{HIKR07}. As a consequence of our methods, we obtain sharp or near sharp results on the distribution of volumes determined by subsets of vector spaces over finite fields and the associated arithmetic expressions. 

In particular, our machinery enables us to prove that if $E \subset {\Bbb F}_q^d$, $d \ge 4$, the $d$-dimensional vector space over a finite field ${\Bbb F}_q$, of size much greater than $q^{\frac{d}{2}}$, and if $E$ is a product set, then the set of volumes of $d$-dimensional parallelepipeds determined by $E$ covers ${\Bbb F}_q$. This result is sharp as can be seen by taking $E$ to equal to $A \times A \times \dots \times A$, where $A$ is a sub-field of ${\Bbb F}_q$ of size $\sqrt{q}$. In three dimensions we establish the same result if $|E| \gtrsim q^{\frac{15}{8}}$. We prove in three dimensions that the set of volumes covers a positive proportion of ${\Bbb F}_q$ if $|E| \ge Cq^{\frac{3}{2}}$. Finally we show that in three dimensions the set of volumes covers a positive proportion of ${\Bbb F}_q$ if $|E| \ge Cq^2$, without any further assumptions on $E$, which is again sharp as taking $E$ to be a $2$-plane through the origin shows. 
\end{abstract} 

\tableofcontents

\section{Introduction}

\vskip.125in

The classical Erd\H os-Szemeredi sum-product problem asks for the
smallest possible value of
$$ \max \{|A+A|, |A \cdot A|\},$$ where $A$ is a finite subset of a
given ring and $|\cdot|$ denotes the cardinality of a finite set.
$$ A+A=\{a+a': a,a' \in A\},$$ and
$$ A \cdot A=\{a \cdot a': a,a' \in A\}.$$

In this case when the ring is real numbers, Erd\H os and Szemeredi
conjectured that
$$ \max \{|A+A|, |A \cdot A|\} \gtrapprox {|A|}^2,$$ where here and throughout the paper, $X \lesssim Y$ means that there exists $C>0$ such that $X \leq CY$. Similarly, $X \lessapprox Y$, with respect to the parameter $N$ means that for every $\epsilon>0$ there exists $C_{\epsilon}>0$ such that $X \leq C_{\epsilon}N^{\epsilon}Y$.  

The best known result, in the setting of real numbers,
$$ \max \{|A+A|, |A \cdot A|\} \gtrapprox {|A|}^{\frac{14}{11}}$$ is
due to Solymosi (\cite{Sol05}), following up on a result due to Elekes
(\cite{El97}).

In finite fields the problem takes on a flavor of its own due to
arithmetic considerations. The first non-trivial result was obtained
by Bourgain, Katz and Tao in \cite{BKT04}. For arbitrary fields the best known results to date are
$$\max \{|A+A|, |A \cdot A|\} \gtrsim \min \{|A|^{\frac{1}{2}}q^{\frac{1}{2}},|A|^2q^{-\frac{1}{2}}\}$$
due to Garaev (\cite{Gr072}). See also \cite{Gr07}, \cite{KS07}, \cite{KS072}, \cite{Sh07}, \cite{Vu07}. 

A related problem is to determine how large $A \subset {\Bbb F}_q$,
the finite field with $q$ elements needs to be in order to assure that
\begin{equation} \label{d2} {\Bbb F}_q^{*} \subset dA^2=A \cdot A+A
\cdot A+\dots+A \cdot A, \end{equation} where ${\Bbb F}_q^{*}$ denotes
the multiplicative group of ${\Bbb F}_q$. Bourgain (\cite{B05}) proved that
(\ref{d2}) holds if $d=3$ and 
$$|A| \ge Cq^{\frac{3}{4}}.$$ See also \cite{BGK06}, \cite{Cr04}, \cite{Gl06}, \cite{GK06}, 
and the references contained therein for related results. 

The second and the third  listed authors developed a combination of geometric and Fourier analytic machinery (\cite{HI07}) to establish the following result. See also \cite{HIKR07}, \cite{HIS07} and \cite{IR07} where this and related machinery is developed in a variety of contexts. 
\begin{theorem} \label{kickass} Let $A \subset {\mathbb F}_q$, where
${\mathbb
F}_q$ is an arbitrary finite field with $q$ elements, such that
$|A|>q^{\frac{1}{2}+\frac{1}{2d}}$. Then
\begin{equation} \label{sex} {\mathbb F}_q^{*} \subset dA^2.\end{equation}

Moreover, suppose that for some constant $C^{\frac{1}{d}}_{size}$,
$$|A| \ge C^{\frac{1}{d}}_{size} q^{\frac{1}{2}+\frac{1}{2(2d-1)}}.$$
Then
\begin{equation} \label{sex2} |dA^2| \ge q \cdot
\frac{C^{2-\frac{1}{d}}_{size}}{C^{2-\frac{1}{d}}_{size}+1}. \end{equation}
\end{theorem}

It follows immediately from Theorem \ref{kickass} that in the case $d=2$,
$$ {\mathbb F}_q^{*} \subset 2A^2 \  \text{if} \ |A|>q^{\frac{3}{4}},$$ and
$$ |2A^2| \ge q \cdot
\frac{C^{\frac{3}{2}}_{size}}{C^{\frac{3}{2}}_{size}+1} \ \text{if} \ 
|A| \ge C_{size}^{\frac{1}{2}} q^{\frac{2}{3}}.$$

This result was proved as a corollary of the following geometric statement. 
\begin{theorem} \label{kickassmama} Let $E \subset {\Bbb F}_q^d$ such that
$|E|>q^{\frac{d+1}{2}}$. Then
$$ {\Bbb F}_q^{*} \subset \{x \cdot y: x,y \in E\}.$$

Suppose that $E \subset {\Bbb F}_q^d$ such that $|E| \ge C_{size}q^{\frac{d}{2}+\frac{d}{2(2d-1)}}$. Suppose, in addition, that $E$ is a product set. Then 
$$ |\{x \cdot y: x,y \in E\}| \ge q\frac{C_{size}^{2-\frac{1}{d}}}{1+C_{size}^{2-\frac{1}{d}}}.$$
\end{theorem}

The main tool used to establish this result is the following geometric incidence theorem. 
\begin{theorem} \label{incidenceprevious} Let $E \subset {\Bbb F}_q^d$
that does not contain the origin. Let
$$ \nu(t)=\sum_{x \cdot y=t} E(x)E(y),$$ where here and throughout the
paper, $E(x)$ denotes the characteristic function of $E$. Then 
$$ \nu(t)={|E|}^2q^{-1}+R(t),$$ where
$$ |R(t)| \leq |E|q^{\frac{d-1}{2}}, \ \text{if} \ t \not=0,$$ and
$$ |R(0)| \leq |E|q^{\frac{d}{2}}.$$

Moreover,
$$ \sum_t \nu^2(t) \leq {|E|}^4q^{-1}$$
$$+|E|q^{2d-1}\sum_k {|\widehat{E}(k)|}^2 |E \cap l_k|,$$ where
$$ l_k=\{tk: t \in {\Bbb F}^{*}_q\}.$$
\end{theorem}

\subsection{Focus of this paper} The purpose of this paper is to develop the geometric incidence machinery to study the distribution of volumes and to apply the resulting estimates to the sum-product type estimates. More precisely, define $vol(x^1, \dots, x^d)$ to be the determinant of the matrix whose rows are $x^j$s. Recall that 
$$vol(x^1, \dots, x^d) =x^1 \cdot (x^2 \wedge x^3 \wedge \dots \wedge x^d),$$ where the dot product is defined by the usual formula
$$ u \cdot v=u_1v_1+u_2v_2+\dots+u_dv_d,$$ and the generalized cross product, sometimes called the wedge product, is given by the identity
$$ u^2 \wedge \dots \wedge u^d=\det \left( \begin{array}{ccc}
i \\
u^2\\
\dots \\
u^d \end{array} \right),$$ where 
$$ i=(i_1, i_2, \dots, i_d),$$ indicating the coordinate directions in ${\Bbb F}_q^d$. Similarly define
$$ vol(E)=\{vol(x^1, \dots, x^d): x^j \in E\}.$$

The question we ask is, under a variety of natural structural assumptions, how large does $E$ need to be to ensure that ${\Bbb F}_q^{*} \subset vol(E)$, or, more modestly, that $vol(E)$ contains a positive proportion of the elements of ${\Bbb F}_q$. 

Taking $E=A \times A \times \dots \times A$, a product set, will allow us to study some special cases of the following general arithmetic problem. Let 
$$ \Lambda: {({\Bbb F}_q^d)}^D={\Bbb F}_q^d \times \dots \times {\Bbb F}_q^d  \to {\Bbb F}_q$$ be a multi-homogeneous form in the sense that if $t_j \in {\Bbb F}_q$, then 
$$ \Lambda(t_1x^1, t_2x^2 \dots, t_Dx^D)=
t_1^{k_1}t_2^{k_2}\dots t_D^{k_D}\Lambda(x^1, \dots, x^D).$$ 

The question we ask is, large does $A \subset {\Bbb F}_q$ need to be so that 
\begin{equation} \label{detset} \Lambda(A)=\{\Lambda(x^1, \dots, x^D): x^j \in A \times \dots \times A\} \end{equation} contains the whole ${\Bbb F}_q$, or at least a positive proportion?

If $D=2$ and 
$$ \Lambda(x^1,x^2)=x^1 \cdot x^2,$$ we are in the realm of Theorem \ref{kickass}. If 
$$ \Lambda(x^1, x^2, \dots, x^d)=x^1_1\dots x^D_1+\dots+x^1_d \dots x^D_d,$$ we would be looking at $DA^d$, the problem recently studied in \cite{GK06}. In this paper we illustrate our method by studying the case 
\begin{equation} \label{determinantdef} \Lambda(x^1, \dots, x^d)=\det(x^1, \dots, x^d), \end{equation} the determinant of the $d$ by $d$ matrix with columns given by $x^j$s, though it will be clear that the method applies to a large variety of multi-homogeneous forms. 

\subsection{Acknowledgements:} The third listed author wishes to thank Boris Bukh, Alexei Glibichuk, Sergei Konyagin, Ilya Shkredov, and Jozsef Solymosi for a number of very useful conversations related to the topics covered by this paper during his visit to the Institute for Advanced Study in November, 2007. 

\vskip.125in 

\section{Statement of results} 

\vskip.125in 

\subsection{Generalized geometric incidence estimates:} 

The main tool in our investigation is the following generalized geometric incidence theorem which can be viewed as a functional version of Theorem \ref{incidenceprevious}. 
\begin{theorem} \label{incidencenow} Let
$$ \nu(t)=\sum_{x \cdot y=t} f(x)g(y),$$ where $f,g$ are non-negative
functions on ${\Bbb F}_q^d$. Then
\begin{equation} \label{pointwiseincidencenow} \nu(t)={||f||}_1{||g||}_1q^{-1}+R(t), \end{equation} where
$$ |R(t)| \leq {||f||}_2{||g||}_2q^{\frac{d-1}{2}} \ \text{if} \ t \not=0.$$

\vskip.125in

Moreover, if $(0, \dots, 0) \notin support(f) \equiv E$, then
\begin{equation} \label{L2incidencenow} \sum_t \nu^2(t) \leq {||f||}_2^2 \cdot |E| \cdot {||g||}_1^2 
\cdot q^{-1} \end{equation}
$$+{||f||}_2^2 \cdot q^{2d-1}
\sum_{k \not=(0, \dots, 0)} {|\widehat{g}(k)|}^2 |E \cap l_k|,$$
where, as before,
$$ l_k=\{tk: t \in {\Bbb F}^{*}_q\}.$$
\end{theorem}

\begin{remark} \label{plusminus}The proof of Theorem \ref{incidencenow} goes through without any essential changes if the dot product $x \cdot y$ is replaced by any non-degenerate bi-linear form 
$B(x,y)$.
\end{remark}

\vskip.125in

\subsection{Distribution of volumes and applications to sums and products} 

Before stating our main results, we need the following definitions. 
\begin{definition} We say that $E \subset {\Bbb F}_q^d$ is {\it product-like} if given any $n$-dimensional sub-space $H_n \subset {\Bbb F}_q^d$, 
$$ |E \cap H_n| \lesssim {|E|}^{\frac{n}{d}}.$$ 
\end{definition} 

\begin{definition} We say that $E \subset {\Bbb F}_q^d$ is in {\it general position} if given any $n$-dimensional sub-space $H_n \subset {\Bbb F}_q^d$, there exist $d-n$ linearly independent vectors of $E$ whose span does not intersect $H_n$. 
\end{definition} 

\begin{remark} The meaning of the general position condition is that if $E \cap H_n$ determines a positive proportion of all the $n$-dimensional volumes, then by elementary geometry, $E$ determines a positive proportion of all the $d$-dimensional volumes. \end{remark} 

Our main results are the following: 

\begin{theorem} \label{mainproductvolume} Suppose that $E \subset {\Bbb F}_q^3$ is product-like. Then 
\begin{equation} \label{mainpointwise} {\Bbb F}_q^{*} \subset vol(E) \ \text{if} \ |E| \ge C
q^{\frac{15}{8}}
\end{equation} with a sufficiently large constant $C>0$. 
\end{theorem} 

\begin{theorem} \label{product3d} Suppose that $E=A \times A \times A$ and $|A|>\sqrt{q}$. Then 
$$ |vol(E)|>\frac{q}{2}.$$ 
\end{theorem} 

\begin{theorem} \label{product4d} Suppose that $E=A \times A \times A \times A$ and $|A|>\sqrt{q}$. Then 
$$ vol(E)={\Bbb F}_q.$$ 
\end{theorem} 

\begin{remark} Observe that Theorem \ref{product3d} and Theorem \ref{product4d} are in general sharp. To see this we can simply take $A$ to be a sub-field. For sharpness examples in the case of prime fields, see \cite{Gl06}. \end{remark} 

\begin{theorem} \label{mainvolume} Suppose that $E \subset {\Bbb F}_q^3$ is in general position and 
$|E| \ge Cq^{2}$ with a sufficiently large constant $C>0$. Then there exists $c>0$ such that 
$$ |vol(E)| \ge cq.$$ 
\end{theorem} 

\begin{remark} The {\it general position} assumption in Theorem \ref{mainvolume} is easily removed since any set $E$ with $|E|>q^2$ is in general position. The assumption $|E| \ge Cq^2$ is sharp since a two-dimensional plane that passes through the origin determines exactly one volume- the zero volume. 
\end{remark}  

\vskip.125in 

\section{Proof of the generalized geometric incidence estimate (Theorem \ref{incidencenow})}

\vskip.125in

We have
$$ \nu(t)=\sum_{x \cdot y=t} f(x)g(y)$$
$$=\sum_{x,y} q^{-1} \sum_s \chi(s(x \cdot y-t)) f(x)g(y)$$
$$={||f||}_1{||g||}_1q^{-1}+q^{-1} \sum_{x,y} \sum_{s \not=0} \chi(s(x
\cdot y-t)) f(x)g(y)$$
$$={||f||}_1{||g||}_1q^{-1}+R(t).$$

Using the Cauchy-Schwartz inequality,
$$ R^2(t) \leq {||f||}_2^2 \cdot q^{-2}
\sum_x \sum_{y,y'} \sum_{s,s' \not=0} g(y)g(y') \chi(x \cdot
(sy-s'y')) \chi(t(s'-s))$$
$$=q^{d-2}{||f||}_2^2 \cdot \sum_{sy=s'y': s,s' \not=0} g(y)g(y') \chi(t(s'-s))$$
$$=q^{d-2}{||f||}_2^2 \cdot \sum_{s \not=0} \sum_y
g^2(y)+q^{d-2}{||f||}_2^2 \cdot
\sum_{sy=s'y'; s \not=s'; s,s' \not=0} g(y)g(y') \chi(t(s'-s))$$
$$=q^{d-2}(q-1){||f||}_2^2{||g||}_2^2+q^{d-2}{||f||}_2^2 \cdot
\sum_{sy=s'y'; s \not=s'; s,s' \not=0} g(y)g(y') \chi(t(s'-s))$$
$$=q^{d-2}(q-1){||f||}_2^2{||g||}_2^2+q^{d-2}{||f||}_2^2 \cdot
\sum_{a \not=0,1;b \not=0} \sum_y g(y)g(ay) \chi(tb(1-a))$$
$$={||f||}_2^2{||g||}_2^2q^{d-2}(q-1)-q^{d-2}{||f||}_2^2 \cdot \sum_{a
\not=0,1} \sum_y g(y)g(ay),$$ and the result follows.

\vskip.25in

To prove the second part of Theorem \ref{incidencenow}, apply
Cauchy-Schwartz once again to see that
$$ \nu^2(t) \leq {||f||}_2^2 \cdot \sum_{x,y,y': x \cdot y=x \cdot
y'=t} E(x)g(y)g(y'),$$ where $E=support(f)$.

It follows that
$$ \sum_t \nu^2(t) \leq {||f||}_2^2 \cdot \sum_{x \cdot y=x \cdot y'}
E(x)g(y)g(y')$$
$$={||f||}_2^2 q^{-1} \sum_{x,y,y'} \sum_s \chi(s(x \cdot y-x \cdot
y')) E(x)g(y)g(y')$$
$$={||f||}_2^2q^{-1}|E|{||g||}_1^2+{||f||}_2^2q^{2d-1} \sum_{s \not=0}
\sum_x {|\widehat{g}(x)|}^2 E(sx)$$
$$={||f||}_2^2q^{-1}|E|{||g||}_1^2+{||f||}_2^2q^{2d-1} \sum_x
{|\widehat{g}(x)|}^2 |E \cap l_x|,$$ as desired.

\vskip.25in

\section{Proof of the volume estimates (Theorem \ref{mainproductvolume} and Theorem \ref{mainvolume})}

\vskip.125in

Let $f(x)=E(x)$ and define
$$ g_0(x)=|\{(u^2, u^3) \in E \times E: u^2 \wedge u^3=x\}|.$$

Observe that
\begin{equation} \label{fnorm} {||f||}_1={||f||}_2^2=|E|. \end{equation}

On the other hand,
\begin{equation} \label{gnorm1} {||g_0||}_1={|E|}^{2}. \end{equation}

Let 
$$ \nu(t)=|\{(x^1, x^2, x^3) \in E \times E \times E: vol(x^1, x^2, x^3)=t\}|$$ and observe that it equals 
$$ \sum_{x \cdot y=t} f(x)g_0(y),$$ where $f$ and $g_0$ are as above. 

\subsection{Elimination of the origin} Let $g(x)=g_0(x)$ if $x \not=(0,0,0)$ and $0$ otherwise. The argument below is simplified considerably if we work with $g$ instead of $g_0$.  We may work with $g$ provided that we show that 
$$ \sum g(x) \gtrsim {|E|}^{2}.$$ To do this, it suffices to show that 
$$ g_0(0,0,0) \leq c{|E|}^{2}$$ for some $0<c<1$. Indeed, 
$$ g_0(0,0,0)=|\{(u^2, u^3) \in E \times E: u^2 \wedge u^3=(0,0,0)\}|$$ 
$$ \leq |E| \cdot \max_{H_1} |E \cap H_{1}|,$$ where the maximum is taken over all the $1$-dimensional sub-spaces of ${\Bbb F}_q^3$. If $E$ is product-like, then 
$$ |E \cap H_1| \lesssim {|E|}^{\frac{1}{3}}$$ and so 
$$ g_0(0,0,0) \lesssim {|E|}^{2-\frac{2}{3}}<<{|E|}^{2}.$$

On the other hand, if 
\begin{equation} \label{malo} |E \cap H_n| \lesssim q^{\frac{n+1}{2}} \end{equation} for every $n$-dimensional sub-space $H_n$, then 
$$ |E \cap H_{1}| \lesssim q,$$ so 
$$ g_0(0,0,0) \leq |E|q,$$ and this quantity is much smaller than ${|E|}^{2}$ if $|E|$ is much larger than $q$. 

Throughout the argument below we shall either prove results about product-like sets, or sets where the condition (\ref{malo}) can be inductively assumed. Thus the origin has indeed been eliminated from the domain of $g_0$ without any harm.  

In the arguments below, we shall work with $g_{H_n}$, defined just like $g$ with respect to sub-spaces $H_n$. We also define $g_H$ to be $0$ at the origin as we do not need to worry about the lower bound in this case. 

\vskip.125in 

\subsection{Proof of Theorem \ref{mainvolume}:} 

\vskip.125in 

We shall need the following estimate, proved in a subsequent section. 
\begin{lemma} \label{gnorm2generic} Suppose that $E \subset {\Bbb F}_q^3$ such that $|E| \gtrsim q^2$ and $|E \cap H_n| \lesssim q^{\frac{n+1}{2}}$ for every $n$-dimensional sub-space $H_n$. Then 
$$ {||g||}_2^2 \lesssim {|E|}^{2} q^2.$$ 
\end{lemma} 

The result holds in two dimensions by Theorem \ref{kickassmama}. We may assume that for any sub-space $H_n$, 
$$ |E \cap H_n| \lesssim q^{\frac{n+1}{2}}$$ for $1 \leq n \leq 2$ for otherwise the induction hypothesis would imply that we recover a positive proportion of all the $n$-dimensional volumes and the definition of general position would then imply that we recover the positive proportion of all the $d$-dimensional volumes. Thus the assumptions of Lemma \ref{gnorm2generic} are satisfied and we see using Cauchy-Schwartz that 
\begin{equation} \label{erectpenis} {|E|}^{6}={\left( \sum_t \nu(t) \right)}^2 \leq |vol(E)| \cdot \sum_t \nu^2(t). \end{equation} 

By Theorem \ref{incidencenow}, keeping in mind that $E$ is not assumed to be product-like, followed by Lemma \ref{gnorm2generic} we see that
$$ \sum_t \nu^2(t) \lesssim |E|q^3 {||g||}_2^2$$
$$ \lesssim |E|q^3 \cdot q^2 {|E|}^{2}.$$ 

Inserting this back into (\ref{erectpenis}) we see that 
$$ |vol(E)| \gtrsim \frac{{|E|}^6}{{|E|}^3 q^5}=\frac{{|E|}^3}{q^5}.$$

This expression is $\gtrsim q$ if 
$$|E| \gtrsim q^2,$$ as desired. This completes the proof of Theorem \ref{mainvolume}. 

\vskip.125in 

\subsection{Proof of Theorem \ref{mainproductvolume}}

We shall need the following estimate proved in a subsequent section. 
\begin{lemma} \label{gnorm2} Suppose that $E \subset {\Bbb F}_q^3$ is product-like. Then 
$$ {||g||}_2^2 \lesssim {|E|}^{\frac{7}{3}}q.$$ 
\end{lemma} 

Applying the estimate (\ref{pointwiseincidencenow}) of Theorem \ref{incidencenow} we see that if $t \not=0$, then 
$$ \nu(t)={|E|}^3q^{-1}+R(t),$$ where 
$$ |R(t)| \leq q {|E|}^{\frac{1}{2}} \cdot {||g||}_2.$$ 

By Lemma \ref{gnorm2},
$$ {|R(t)|}^2 \lesssim q^{3} {|E|}^{3+\frac{1}{3}},$$ provided that 
$$ |E| \lesssim q^{2}.$$ 

It follows that $\nu(t)>0$ if 
$$ {|E|}^{\frac{8}{3}} \ge Cq^5$$ with a sufficiently large constant $C>0$. It follows that $\nu(t)>0$ if 
$$ |E| \ge Cq^{\frac{15}{8}},$$ which completes the proof of the point-wise estimate. 

\vskip.125in 

\section{Proof of the key estimate (Lemma \ref{gnorm2})} 

\vskip.125in 

We shall make use of the following calculations. 
\begin{lemma} \label{dveruki} Let $E \subset {\Bbb F}_q^3$ be product-like. Let $H_2$ be a $2$-dimensional subspace of ${\Bbb F}_q^3$.  Then 
\begin{equation} \label{planeaverage} \sum_{H_2} {|E \cap H_{2}|}^{2} \lesssim {|E|}^2 \end{equation} provided that 
$$ |E| \gtrsim q^{\frac{3}{2}}.$$ 
\end{lemma} 

\begin{lemma} \label{dverukigeneric} Let $E \subset {\Bbb F}_q^3$. Let $H_2$ be an $2$-dimensional sub-space of ${\Bbb F}_q^3$. Then 
\begin{equation} \label{planeaveragegeneric} \sum_{H_2} {|E \cap H_{2}|}^{2} \lesssim {|E|}^2 \end{equation} provided that 
$$ |E| \gtrsim q^{2}.$$ 
\end{lemma} 

\subsubsection{Proof of Lemma \ref{dveruki} and Lemma \ref{dverukigeneric}:} To prove the estimate (\ref{planeaverage}), observe that 
$$ \sum_{H_2} {|E \cap H_2|}^{2}={(q-1)}^{-1} \sum_{x} \left( \sum_{x \cdot y=0} E(y) \right)^{2}$$
$$={(q-1)}^{-1} \sum_{x \in {\Bbb F}_q^3} \Pi_{j=1}^{2} \sum_{x \cdot y^j=0} E(y^j).$$

If $E$ is product-like, then this expression is
\begin{equation} \label{turdproduct} \lesssim {|E|}^2+q{|E|}^{\frac{4}{3}} \lesssim {|E|}^2 \end{equation} if $|E| \ge Cq^{\frac{3}{2}}$, as claimed. 

Under the assumptions of Lemma \ref{dverukigeneric}, this expression is 
\begin{equation} \label{turdgeneric} \lesssim {|E|}^{2}+q^2 \cdot |E| \lesssim {|E|}^2 \end{equation} if $|E| \ge Cq^2$, as desired. 

In both (\ref{turdproduct}) and (\ref{turdgeneric}) we use the fact that if $y^1, y^{2}$ span the $1$-dimensional sub-space, 
$$ |\{x: x \cdot y^j=0; \ j=1,2 \}| \approx q^2.$$ 

In order to establish (\ref{turdproduct}) we also use the fact that the restriction of $E$ to a $1$-dimensional subspace has size $\lesssim {|E|}^{\frac{1}{3}}$, whereas in the general case we simply use the fact that the size of this intersection does not exceed $q$. This completes the proof of the estimate (\ref{planeaverage}) and (\ref{planeaveragegeneric}). 
\vskip.125in 

\subsection{The conclusion of the proof of Lemma \ref{gnorm2} and Lemma \ref{gnorm2generic}:} 

We have
$$ {||g||}_2^2 \lesssim \sum_j \sum_{H_{2} \in G_0(2,3)} \nu^2_{H_{2}}(j)$$
\begin{equation} \label{basicreduction} \lesssim \sum_{H_{2}} {|E \cap H_{2}|}^{4}q^{-1}+II, \end{equation} where 
$$ I=\sum_{H_{2}} {|E \cap H_{2}|}^{4}q^{-1}$$ and 
\begin{equation} \label{IIproduct} II={II}_{prod} \lesssim {|E|}^{\frac{1}{3}} \sum_{H_{2}}  |E \cap H_{2}| 
q {||g_{H_{2}}||}_{L^2(H_{2})}^2 \end{equation} if $E$ is product-like and 
\begin{equation} \label{IIgeneric} II={II}_{generic} \lesssim q\sum_{H_{2}}  |E \cap H_{2}| q {||g_{H_{2}}||}_{L^2(H_{2})}^2 \end{equation} if $E$ is in general position. 

\subsection{Estimation of $I$:} Suppose that $E$ is product like and $|E| \gtrsim q^{\frac{3}{2}}$. Since 
$$ {|E \cap H_{2}|}^{2} \lesssim {|E|}^{\frac{4}{3}}$$ due to the fact that $E$ is a product set and $H$ is a $2$-dimensional plane, Lemma \ref{dveruki} gives us 
$$ I \lesssim q^{-1} {|E|}^{\frac{4}{3}} \cdot {|E|}^{2}$$
$$ \lesssim {|E|}^{\frac{7}{3}} q$$ as long as 
$$ |E| \lesssim q^2$$ as desired. 

If $E$ is in general position, 
$$ {|E \cap H_{2}|}^{2} \lesssim q^3,$$ so 
$$ I \lesssim q^{-1}q^3 {|E|}^{2}$$ by Lemma \ref{dverukigeneric}.

\vskip.125in

\subsection{Estimation of $II$:} Since $g_{H_2}$ is simply the characteristic function of $E \cap H_2$, Lemma \ref{dveruki} implies that 
$$ {||}_{prod} \lesssim q{|E|}^{\frac{7}{3}},$$ if $|E| \gtrsim q^{\frac{3}{2}}$, as desired. 

Similarly, Lemma \ref{dverukigeneric} implies that 
$$ {II}_{generic} \lesssim q^2{|E|}^2,$$ if $|E| \gtrsim q^2$ and we are done. 

\vskip.125in 

\section{Proof of Theorem \ref{product3d} and \ref{product4d}} 

\vskip.125in 

We shall the following lemmae. See \cite{BKT04} and \cite{Gl06} 
\begin{lemma}\label{Glib} Suppose, $A\subset {\mathbb F}_q^*$ is such
that $|A|^2>q$. Then there exist elements $\alpha,\beta\in A-A$, such
that \begin{equation}\label{Gl}|\alpha A\pm \beta A|>\frac{q}{2}.\end{equation}\end{lemma}

\begin{lemma} \label{two} If $C\subset{\mathbb F}_q$ is such that
$|C|>\frac{q}{2}$, then $C\pm C={\mathbb F}_q$.
\end{lemma}

Using Lemma \ref{Glib}, we have the following.
\begin{lemma}\label{kl}
Suppose, $A\subset {\mathbb F}_q^*$ is such that $|A|^2>q$. Let $B=A-A$.
Then $B^2 - B^2={\mathbb F}_q.$\end{lemma}
{\sf Proof.} The claim is that the set $D$ of two by two determinants
$$\left|
  \begin{array}{cc}
    x_1 & x_2 \\
    y_1 & y_2 \\
  \end{array}
 \right|$$
 with elements in $B$ cover ${\mathbb F}_q$. Let us fix
$x_1=\alpha,\,x_2=\beta$, from Lemma \ref{Glib}, which are determined by
$A$ only. Let $y_1=u_1-v_1,\,y_2=u_2-v_2$, where $u_1,v_1,u_2,v_2\in A$. Let
 $C=\alpha A - \beta A$. We have
 $$
 D= C-C,
 $$
 and the result follows from the fact that $|C|>\frac{q}{2}$ by
Lemma \ref{Glib} and Lemma \ref{two}.

We are now ready to prove Theorem \ref{product4d}. Consider the determinants in the form

 $$\left|
    \begin{array}{cccc}
      x_1 & x_2 & x_3 & x_4 \\
      y_1 & y_2 & y_3 & y_4 \\
      u_1 & u_2 & x_3 & x_4 \\
      v_1 & v_2 & y_3 & y_4 \\
    \end{array}
  \right| =
  \left|
    \begin{array}{cccc}
      x_1-u_1 & x_2-u_2 & 0 & 0 \\
      y_1-v_1 & y_2 -v_2 & 0 & 0 \\
      u_1 & u_2 & x_3 & x_4 \\
      v_1 & v_2 & y_3 & y_4 \\
    \end{array}
  \right|$$
  $$
   =  (x_3y_4-y_3x_4) \left|
    \begin{array}{cccc}
      x_1-u_1 & x_2-u_2  \\
      y_1-v_1 & y_2 -v_2 \end{array}
  \right|.
 $$
 The statement now follows from Lemma \ref{kl}.

We now prove Theorem \ref{product3d}. Consider the determinants in the form $$\left|
    \begin{array}{cccc}
      x_1 & x_2 & x_3  \\
      y_1 & y_2 & y_3  \\
      u_1 & u_2 & x_3
    \end{array}
  \right| =
  \left|
    \begin{array}{cccc}
      x_1-u_1 & x_2-u_2 & 0  \\
      y_1 & y_2 &y_3\\
      u_1 & u_2 & u_3
    \end{array}
  \right|$$
  $$
   =  (x_1-u_1)(y_2u_3 - y_3u_2) - (x_2-u_2) (y_1u_3 - y_3u_1)
 $$
 $$
 = u_3 [ (x_1-u_1)y_2 - (x_2-u_2)y_1] - y_3[(x_1-u_1)u_2 - (x_2-u_2)u_1]
 $$
Let  $(x_1-u_1)=\alpha,\,(x_2-u_2)=\beta$ come from Lemma \ref{Glib},
having therefore fixed
 $x_1,u_1,x_2,u_2$. Now fix $u_3\neq0$, and some $y_3$.
 The statement now follows from Lemma \ref{Glib}.

\enddocument